\documentclass[12pt]{article}
\title{Uniqueness of Crepant Resolutions and Symplectic Singularities}
\author{Baohua Fu and Yoshinori Namikawa}
\usepackage{amsmath,amssymb,amsthm, amscd}
\chardef\bslash=`\\
\newtheorem{Thm}{Theorem}[section]
\newtheorem{Cor}[Thm]{Corollary}
\newtheorem{Lem}[Thm]{Lemma}
\newtheorem{Prop}[Thm]{Proposition}
\newtheorem{Def}{Definition}
\newtheorem{Rque}[Thm]{Remark}
\newtheorem{Conj}{Conjecture}
\newtheorem{Exam}[Thm]{Example}

\def\cit{{\mathbb C}}

\def\qit{{\mathbb Q}}
\def\pit{{\mathbb P}}

\def\0{{\mathcal O}}
\def\g{{\mathfrak g}}

\begin{document}
\maketitle
{\em Running title:} Uniqueness of crepant resolutions.\\

{\em Abstract:}
We prove the uniqueness of crepant resolutions for some quotient singularities and for some nilpotent orbits. 
The finiteness of non-isomorphic symplectic resolutions for 4-dimensional symplectic singularities is
proved. We also give an example of symplectic singularity which admits two non-equivalent symplectic resolutions.\\


{\em Classification:} 14E15 \\

{\em Titre en fran\c{c}ais:} Unicit\'e des r\'esolutions cr\'epantes et singularit\'es symplectiques. \\

{\em R\'esum\'e:} Nous d\'emontrons l'unicit\'e des r\'esolutions cr\'epantes pour certaines singularit\'es quotient
et pour certaines adh\'erences d'orbites nilpotentes. La finitude des r\'esolutions symplectiques non-isomorphes
pour les singularit\'es symplectiques de dimension 4 est d\'emontr\'ee. Nous construisons aussi un exemple de
singularit\'e symplectique qui admet deux r\'esolutions symplectiques non-\'equivalentes. \\

{\em Mots cl\'es:} r\'esolutions cr\'epantes,  singularit\'es symplectiques\\

{\em Classification:} 14E15

\section{Introduction}
In this paper, we work over the field  $\cit$ of complex numbers.
Let $W$ be an algebraic variety, smooth in codimension 1, such that $K_W$ is a Cartier divisor.
 Recall that a resolution of singularities $\pi: X \to W$ is called {\em crepant}
if $\pi^* K_W = K_X.$ In this note, we will only consider projective crepant resolutions, i.e. $\pi$ is projective.
   Let $\pi^+: X^+ \to W$ be another (projective) crepant resolution of $W$. 
\begin{Def}
(i) $\pi$ and $\pi^+$ are said {\em isomorphic} if the natural birational map $\pi^{-1} \circ \pi^+: X^+ --\to X $ is an isomorphism;

(ii) $\pi$ and $\pi^+$ are said {\em equivalent} if there exists an automorphism $\psi$ of $W$ such that $\psi \circ \pi$ and 
$\pi^+$ are isomorphic.
\end{Def}

As easily seen, any two crepant resolutions of A-D-E singularities are isomorphic.
The purpose of this note is to study projective crepant resolutions (mostly for symlectic singularities) up to isomorphisms and up to equivalences.

A special case of crepant resolutions is symplectic resolutions for symplectic singularities. Following \cite{Bea}, a variety $W$, 
smooth in codimension 1,  
is said to have {\em symplectic singularities} if there 
exists a holomorphic
symplectic 2-form $\omega$ on $W_{\mathrm{reg}}$ such that for any resolution of singularities $\pi: X \rightarrow W$, the 2-form
$\pi^* \omega$ defined a priori on $\pi^{-1}(W_{\mathrm{reg}})$ can be extended to a holomorphic 2-form on $X$. If furthermore 
the 2-form $\pi^* \omega$ extends to a holomorphic symplectic 2-form on the whole of $X$ for some resolution of $W$, 
then we say that $W$ admits a {\em symplectic resolution}, and the resolution $\pi$ is called {\em symplectic}.

For a symplectic singularity, a resolution is symplectic if and only if it is crepant (see for example Proposition 1.1 \cite{Fu}).
In recent years, there appeared many studies on symplectic resolutions for symplectic singularities (see \cite{CMS}, \cite{Fu}, \cite{Fuj},
\cite{Ka}, \cite{Kal}, \cite{Nam} and \cite{Wie} etc.).
 
Our first theorem on uniqueness of crepant resolutions is the following: \vspace{0.2cm}

{\bf Theorem (2.2)}. Let $W_i, i=1, \cdots, k$ be normal locally $\qit$-factorial singular varieties which admit
 a crepant resolution $\pi_i: X_i \to W_i$ such that 
$E_i: = Exc(\pi_i)$ is an irreducible divisor. Suppose that $W: = W_1 \times \cdots \times W_k$ is locally $\qit$-factorial.
Then any crepant resolution of $W$ is isomorphic to the product 
$$\pi:= \pi_1 \times \cdots \times \pi_k:  X:= X_1 \times \cdots \times X_k \to W_1  \times \cdots \times W_k.$$ 

It applies to many varieties with quotient singularities. For example it shows that for any smooth surface $S$, its $n$th symmetric product 
$S^{(n)}$ admits a unique crepant resolution, which is  given by the Douady-Barlet resolution:  $S^{[n]} \to S^{(n)}.$
As to the nilpotent orbit closures, we have \vspace{0.2cm}

{\bf Theorems (3.1)}. Let $\0$ be a nilpotent orbit in a complex semi-simple Lie algebra $\g$. Then $\overline{\0}$ admits at most finitely
many non-isomorphic symplectic resolutions.

This result is an easy corollary of our previous work in \cite{Fu}.
Some other partial results are also presented in Section 3. The above theorem motivates the following:\vspace{0.2cm}

{\bf Conjecture (1)}. Any symplectic singularity admits at most finitely many non-isomorphic symplectic resolutions.

In Section 4, we prove this conjecture in the 4-dimensional case. As to the relation between two symplectic resolutions, we have the following:
\vspace{0.2cm}

{\bf Conjecture (2)}. Let $W$ be a normal symplectic singularity. Then for any two symplectic resolutions $f_i: X_i \to W, i=1, 2$, 
there are deformations $\mathcal{X}_i \stackrel{F_i}\to 
\mathcal{W}$ of $f_i$ such that, for $s \in S\setminus 0$, $F_{i,s}: 
\mathcal{X}_{i,s} \to \mathcal{W}_s$ are isomorphisms. In particular, 
$X_1$ and $X_2$ are deformation equivalent. 

 By constructing explicitly the deformations,
we prove this conjecture for symplectic resolutions of nilpotent orbit closures
in $\mathfrak{sl}(n)$ in section 4.

 Finally in Section 5, we construct an example of symplectic singularity of dimension 4
which admits two non-equivalent symplectic resolutions. 

The following proposition gives some applications of results presented in this note. 
\begin{Prop}\label{Prop} 
Let $W$ be an algebraic variety, smooth in codimension 1.
If up to isomorphisms, $W$ admits a unique crepant resolution $\pi: X \to W$, then any automorphism of $W$ lifts to $X$.
\end{Prop}
\begin{proof}
 Let $\psi: W \to W$ be an automorphism. Then $\psi \circ \pi: X \to W$ is again a crepant resolution, which is isomorphic to $\pi$ by 
hypothesis, thus there exists an automorphism $\widetilde{\psi}$ of $X$ lifting $\psi$.
\end{proof}

{\em Acknowledgements.} The first named author wants to thank A. Beauville, M. Brion,  A. Hirschowitz and  D. Kaledin for helpful
 discussions and suggestions.
 
\section{Quotient singularities}
\begin{Lem} \label{Lem1}
Let $W$ be a normal locally $\qit$-factorial variety and $\pi: X \to W$ a projective resolution. 
Then $\mathrm{Exc}(\pi)$ is of pure codimension $1$ and 
if $\mathrm{Exc}(\pi) = \cup_{i=1}^n E_i$ is the decomposition into irreducible components,
 then $\0_X(-\sum_i a_i E_i)$ is $\pi$-ample for some $a_i > 0$. 
\end{Lem}
\begin{proof}
The first claim is well-known (see 1.40 \cite{Deb}), which follows from the normality and $\qit$-factority of $W$. 
For the second claim, by 1.42 \cite{Deb}, $\0_X(-\sum_i a_i E_i)$ is $\pi$-very ample for some $a_i \geq 0$.
Suppose that $a_{i_0} = 0$, then take a point $x \in E_{i_0} - \cup_{i \neq i_0} E_i$ and a $\pi$-exceptional
curve $C$ passing $x$. Note that $C$ is not contained in $E_i, i\neq i_0$, thus $C \cdot E_i \geq 0, i\neq i_0.$
This gives $ (-\sum_i a_i E_i) \cdot C \leq 0,$ which is absurd since  $\0_X(-\sum_i a_i E_i)$ is $\pi$-very ample.
\end{proof}
We are indebted to M. Brion for pointing out the reference \cite{Deb}.
\begin{Thm}\label{Thm}
Let $W_i, i=1, \cdots, k$ be normal locally $\qit$-factorial singular varieties which admit a crepant resolution $\pi_i: X_i \to W_i$ such that 
$E_i: = \mathrm{Exc}(\pi_i)$ is an irreducible divisor. Suppose that $W: = W_1 \times \cdots \times W_k$ is locally $\qit$-factorial.
Then any crepant resolution of $W$ is isomorphic to the product 
$$\pi:= \pi_1 \times \cdots \times \pi_k:  X:= X_1 \times \cdots \times X_k \to W_1  \times \cdots \times W_k.$$ 
\end{Thm}
\begin{proof}
The $\pi$-exceptional locus consists of $k$ irreducible divisors $F_i: = X_1 \times \cdots \times E_i \times \cdots \times X_k.  $
We first prove that $-F_i$ is $\pi$-nef for all $i$. Let $C$ be a curve in $X$ such that $\pi(C)$ is a point. Consider the following composite
$$X = X_1 \times \cdots \times X_k \xrightarrow{p_i} X_i \xrightarrow{\pi_i} W_i.$$

Note that $F_i = p_i^*(E_i)$. If $p_i(C)$ is a point $Q$, then $(C,F_i) =  0.$ 
If $p_i(C)$ is a curve, then $$(C, F_i) = [C: p_i(C)] (p_i(C), E_i).$$ Applying Lemma \ref{Lem1} to the resolution 
$\pi_i: X_i \to W_i$, we see that $-aE_i$ is $\pi_i$-ample for some $a>0$, thus $ (p_i(C), E_i)< 0$, since $\pi_i(p_i(C))$ is a point.
Therefore, $-F_i$ is $\pi$-nef.

Assume now that there is another crepant resolution $\pi^+: X^+ \to W.$  Then $X$ and $X^+$ are isomorphic in codimension 1 because $\pi$ and 
$\pi^+$ are both crepant resolutions. In particular, $\mathrm{Exc}(\pi^+)$ contains exactly $k$ irreducible divisors, say $F_i^+, 1 \leq i \leq k.$
Now apply Lemma \ref{Lem1}, $L^+: = \0_{X^+}(- \sum_i a_i F_i^+) $ is $\pi^+$-ample for some $a_i > 0$. Its proper transform by the 
birational map $X^+ --\to X$  coincides with  $L: = \0_{X}(- \sum_i a_i F_i), $ which is $\pi$-nef.

Since $L$ is $\pi$-nef, $\pi$-big and $\pi$ is crepant, the Base Point Free theorem 
 implies that $L^{\otimes m}$ is $\pi$-free for a sufficiently large $m$.  So there is a birational morphism 
$X \to \mathrm{Proj}_W(\oplus_k \pi_*L^{\otimes mk}). $  On the other hand, since $X$ and $X^+$ are isomorphic in codimension 1, 
there is an isomorphism $\pi_*L^{\otimes mk} \simeq \pi^+_*{L^+}^{\otimes mk}. $ Therefore we have a birational morphism 
$X \to X^+$ over $W$. Since $X$ and $X^+$ are both crepant resolutions of $W$, this birational morphism should be an isomorphism over
$W$. Hence $\pi$ and $\pi^+$ are isomorphic.
\end{proof}

For a smooth surface $S$, we denote by $S^{(n)}$ its 
symmetric n-th products (the Barlet space parametrizing 
$0$ cycles on $S$ of length $n$), and we denote by 
$S^{[n]}$ the Douady space parametrizing $0$-dimensional 
subspaces of $S$ with length $n$. 

\begin{Cor}
Let $S_i, i=1, \cdots k$ be a smooth surface. Then any crepant resolution of $S_1^{(n_1)} \times \cdots \times S_k^{(n_k)}$
is isomorphic to the Douady-Barlet resolution
$$ S_1^{[n_1]} \times \cdots \times S_k^{[n_k]} \to     S_1^{(n_1)} \times \cdots \times S_k^{(n_k)}             .$$
\end{Cor}
\begin{Cor}\label{Cor2}
Let $V$ be a symplectic vector space and $G$ a finite subgroup of $Sp(V)$. Suppose that the symplectic reflections of $G$ (i.e. $g \in G$
such that $Fix(g)$ is of codimension 2) form a single conjugacy class. Then any two crepant resolutions of $V/G$ are isomorphic.
\end{Cor}
\begin{proof}
Let $\pi: X \to V/G$ be a crepant resolution. By McKay correspondence proved by D. Kaledin (\cite{Ka3}), there is a one-to-one
correspondence between the conjugacy classes of symplectic reflections in $G$ and closed  irreducible sub-varieties $E$  of codimension 1 in  $X$ such that
$\mathrm{codim}(\pi(E)) = 2$. Notice that such $E$ is exactly irreducible components of  $\mathrm{Exc}(\pi)$.
 By the hypothesis, there is only one such conjugacy class, thus $\mathrm{Exc}(\pi)$ is irreducible.  
\end{proof}
Combining this corollary with Proposition \ref{Prop}, we have immediately the following:
\begin{Cor}\label{Cor3}
Let $V$ be a symplectic vector space and $G$ a finite subgroup of $Sp(V)$. Suppose that the symplectic reflections of $G$ 
form a single conjugacy class and $\pi: X \to V/G$ is a crepant resolution. Then any action of an algebraic group $H$ on $V/G$
 lifts to an $H$-action on $X$.
\end{Cor}
\begin{Rque}
Corollary \ref{Cor2}  gives a generalization of a result proved by D. Kaledin (Theorem 1.9 \cite{Ka}) and Corollary \ref{Cor3}
strengthens Theorem 1.3 {\em loc. cit.}.
\end{Rque}

\begin{Exam} \label{exam}
 Here is an example to show  the condition in Corollary \ref{Cor2} that symplectic  reflections of $G$ form a single conjugacy class
is necessary.  This example has also been considered by A. Fujiki (\cite{Fuj}).

Let $(x,y,z,w)$ be the coordinates of $\cit^4$.  Let $G$ be the subgroup of 
$\mathrm{Aut}(\cit^4)$  generated by three elements  $$ \sigma_1: (x,y,z,w) \to (x,y,-z,-w),$$ 
$$ \sigma_2: (x,y,z,w) \to (-x,-y,z,w),$$  $$ \tau : (x,y,z,w) \to (z,w,x,y). $$ 
Then $G$ is dihedral group of order 8.  Since all elements of $G$ preserves 
the two form $dx \wedge dy + dz \wedge dw$,  the quotient $W := \cit^4/G$ is a  symplectic singularity. 

One sees easily that $W = \mathrm{Sym}^2(\bar{S})$, where $\bar{S} = \cit^2 / \pm 1.$ 
Let $S \to \bar{S}$ be the minimal resolution. 
Let $C$ be its exceptional curve. $C \cong 
\pit^1$ and $(C^2)_S = -2$. Now we have 
a sequence of birational maps 
$$ X := \mathrm{Hilb}^2(S) 
\stackrel{f_1}\to \mathrm{Sym}^2(S) 
\stackrel{f_2}\to \mathrm{Sym}^2(\bar{S}) = W. $$ 
Let $f: X \to W$ be the composite of the maps, which 
is a symplectic resolution of $W$. 
Note that $f_2^{-1}(0) = \mathrm{Sym}^2(C)(\cong 
\pit^2)$. Let $\Delta_C \subset 
\mathrm{Sym}^2(C)$ be the diagonal. Put $F := 
f_1^{-1}(\Delta_C)$. Then $F$ is a $\pit^1$ 
bundle over $\Delta_C (\cong \pit^1)$. 
It can be checked that $F$ is isomorphic to 
the Hirzebruch surface $\Sigma_4$. 
As a consequence, we have $$ f^{-1}(0) = \pit^2 \cup F, $$ 
where $\pit^2$ is the proper transform 
of $\mathrm{Sym}^2(C)$ by $f_1$.  The intersection 
$\pit^2 \cap F$ is a conic of $\pit^2$ 
and, at the same time, is a negative section of 
$F \cong \Sigma_4$.   

$Sing(W)$ has two components $T_1$ (diagonal of $\bar{S} \times \bar{S}$) and $T_2$. 
Let $E_i = f^{-1}(T_i), i=1,2$. Then $f^{-1}(0) \subset E_1, f^{-1}(0) \cap E_2 = F.$
In particular, we see that the resolution $f$ is not symmetric with respect to $T_1$ and $T_2$.

Consider the map $u: \cit^4 \to \cit^4$ defined by $u(x,y,z,w) = (x-z, y-w, x+z, y+w).$
One verifies that $$ u \circ \sigma_1 = \tau \circ u; \quad u \circ \tau  = \sigma_2 \circ u; \qquad u \circ \sigma_2 = \sigma_1 \circ \sigma_2 \circ \tau
\circ u.    $$ Thus $u$ gives an automorphism $\bar{u}$ on $W = \cit^4/G,$ which interchanges $T_1$ and $T_2$.
So the two crepant resolutions $f$ and $f':=\bar{u} \circ f$ are not isomorphic, though they are equivalent.

In fact one can show that the birational map: $(f')^{-1} \circ f:  \mathrm{Hilb}^2(S) --\to \mathrm{Hilb}^2(S)$ is exactly the Mukai flop along
the subvariety $\pit^2 $ of $ \mathrm{Hilb}^2(S).  $
\end{Exam}
 
\section{Nilpotent orbits}
Let $\g$ be a semi-simple complex Lie algebra and $\0$ a nilpotent orbit in $\g$. Then $\overline{\0}$ is singular and smooth in codimension 1.
Let $\widetilde{\0}$ be its normalization, which is a normal variety with  symplectic singularities (\cite{Bea}).
It is proved in \cite{Fu} that any projective  symplectic resolution of $\widetilde{\0}$ is isomorphic to the collapsing of the zero section of 
$T^*(G/P)$ for some parabolic subgroup $P$ of $G$, where $G$ is the adjoint group of $\g$. Notice that $G$ has only finitely many conjugacy classes of 
parabolic subgroups, thus we get
\begin{Thm}\label{nilthm}
Let $\g$ be a complex semi-simple Lie algebra and $\0$ a nilpotent orbit in $\g$. Then $\widetilde{\0}$ admits at most finitely many symplectic
 resolutions, up to isomorphisms. 
\end{Thm}
Notice that any two Borel subgroups in a semi-simple Lie group are conjugate, thus we have 
\begin{Cor}
Let $\mathcal{N}$ be the nilpotent cone of a semi-simple complex Lie algebra $\g$. Then any symplectic resolution of
 $\mathcal{N}$ is isomorphic to the Springer resolution $T^*(G/B) \to  \mathcal{N},$ where $B$ is a Borel 
subgroup of $G$.
\end{Cor}

As to the uniqueness up to isomorphisms of symplectic resolutions
for a nilpotent orbit closure, we have following partial results.
\begin{Prop}
Let $\g$ be a simple complex Lie algebra not of type $A$ and $\0$ a nilpotent orbit in $\g$.
Suppose that $\overline{\0} - \0 = \overline{C}$ for some nilpotent orbit $C$ of codimension 2 in $\0$.
If the singularity $(\overline{\0}, C)$ is of type $A_1$, then any two
symplectic resolutions for $\overline{\0}$ are isomorphic.
\end{Prop}
\begin{proof}
Let $\pi: X \to \overline{\0}$ be a symplectic resolution, then over $U:=C \cup \0$, $\pi$ is isomorphic to the blow-up
of $U$ at $C$, since  $(\overline{\0}, C)$ is of type $A_1$. By the semi-smallness of symplectic resolutions (Proposition 1.4 \cite{Nam} or 
Proposition 1.2 \cite{Ka3}), 
$\mathrm{codim}(\pi^{-1}(\overline{\0} - U)) \geq 2$, thus 
 $\mathrm{Exc}(\pi)$ consists of one irreducible
divisor.  Since $\g$ is not of type $A_k$, $\overline{\0}$ is
$\qit$-factorial (\cite{Fu}). Moreover, the $\pi$-exceptional fiber over $C$ is isomorphic to $\pit^1$, thus connected, so
$\overline{\0}$ is normal (Theorem 1. \cite{KP}). 
 Now the proposition follows from Theorem \ref{Thm}.
\end{proof}

Then one can use results of H. Kraft and C. Procesi in \cite{KP} to determine all nilpotent orbits
which satisfy the hypothesis of the proposition above.
For example, in $\mathfrak{so}(5)$, we find $\0_{[3,1,1]}$. In $\mathfrak{sp}(4)$, we have
$\0_{[2,2]}$. In $\mathfrak{so}(8)$, we
get $\0_{[3,3,1,1]}$, $\0_{[3,1^5]}$ and  $\0_{[2,2,2,2]}$ etc..
\begin{Prop}
Let $\0$ be a nilpotent orbit in $\mathfrak{sl}(n+1, \cit)$. Let {\bf d}= $[d_1,\cdots, d_s]$ be its Jordan decomposition type.
If $d_1=\cdots=d_s,$  then up to isomorphisms, $\overline{\0}$ admits a unique symplectic resolution.
\end{Prop}
\begin{proof}
It is well-known that the closure of any nilpotent orbit  in $\mathfrak{sl}(n+1, \cit)$ is normal and admitting a symplectic resolution.
If $d_1=\cdots=d_s,$ then all polarizations of $\0$ (i.e. parabolics $P$ such that $T^*(G/P)$ is birational to  $\overline{\0} $)    
form a single conjugacy class (see for example Theorem 3.3 (b) \cite{Hes}), thus $\overline{\0}$ admits a unique symplectic resolution, up to isomorphisms.
\end{proof}
\begin{Prop}
Let $\0$ be a nilpotent orbit in a complex  simple Lie algebra of type $B-C-D$, with Jordan decomposition type  {\bf d}= $[d_1,\cdots, d_s]$.
Suppose that:

(i). either there exists some integer $k \geq 1 $ such that $ d_1=\cdots=d_s = 2k$;

(ii). or there exist some integers $q \geq 1, k \geq 1$ such that $ d_1=\cdots=d_q = 2k+1$ and 
$d_{q+1} = \cdots = d_s = 2k$. 

 Then $\overline{\0}$ admits a unique symplectic resolution, up to isomorphisms.
\end{Prop}
\begin{proof}
By Proposition 3.21 and Proposition 3.22 \cite{Fu}, such a nilpotent orbit $\overline{\0}$ admits a symplectic resolution.
Furthermore by the proofs there (see also \cite{Fu2}), two polarizations of $\0$ have conjugate Levi factors. Thus
the number of conjugacy classes of polarizations is given by $N_0$ of Theorem 7.1 (d) \cite{Hes}, which equals to 1 in our case.
Thus $\overline{\0}$ admits a unique symplectic resolution, up to isomorphisms.
\end{proof}

Now we study  symplectic resolutions up to equivalences for a nilpotent orbit $\0:= \0_{[d_1, \cdots, d_k]}$ contained in $\mathfrak{sl}(n)$, where
 $[d_1, \cdots, d_k]$ is the Jordan decomposition type of $\0$. Let $[s_1, \cdots, s_m]$ be the dual partition of $[d_1, \cdots, d_k]$.
We denote by $P_{i_1, ..., i_k} \subset SL(n)$ a parabolic subgroup of flag type  $(i_1, ..., i_k)$; that is, $P_{i_1, ..., i_k}$ stabilizes 
a flag $0 = V_0 \subset V_1 \subset ... \subset V_k = \mathbf{C}^n$ such that $\dim V_j/V_{j-1} = i_j$ for $1 \leq j \leq k$. 
This is equivalent to saying that  $SL(n)/P_{i_1, .... i_k}$ is the flag manifold 
$F(n, n - i_k, n - i_k - i_{k-1}, ..., i_1)$. 

It is well-known that all polarizations of $\0$ are of the form $P_{s_{\sigma(1)}, \cdots, s_{\sigma(m)}}$ for some permutation 
$\sigma \in \Sigma_m$ (see for example Theorem 3.3 \cite{Hes}).
\begin{Prop}\label{sl}
The two symplectic resolutions $$T^*(SL(n)/P_{s_{\sigma(1)}, \cdots, s_{\sigma(m)}}) \to \overline{\0},  \quad 
T^*(SL(n)/P_{s_{\sigma(m)}, \cdots, s_{\sigma(1)}}) \to \overline{\0} $$    
are equivalent.
\end{Prop}
\begin{proof}
Take the dual flags, we get an isomorphism between $ SL(n)/P_{s_{\sigma(1)}, \cdots, s_{\sigma(m)}}$ and $SL(n)/P_{s_{\sigma(m)}, \cdots, s_{\sigma(1)}}$.
Furthermore $\overline{\0}$ is normal. Now the proposition follows from the following lemma.
\end{proof}
\begin{Lem} \label{equiv}
Let $W$ be an affine normal variety and $\pi_i: X_i \to W, i=1,2,$ two crepant resolutions. Then $\pi_1$ is equivalent to
$\pi_2$ if and only if $X_1$ is isomorphic to $X_2$.
\end{Lem}
\begin{proof}
The isomorphism $X_1 \cong X_2$  induces an isomorphism of $\cit$-algebras 
$\Gamma (X_1, \mathcal{O}_{X_1}) 
\cong \Gamma (X_2, \mathcal{O}_{X_2}),$ 
thus an isomorphism of algebraic varieties $$\mathrm{Spec}(\Gamma (X_1, \mathcal{O}_{X_1}))
\cong \mathrm{Spec}(\Gamma (X_2, \mathcal{O}_{X_2})).$$ 
The morphism $\pi_i$ gives an injective morphism from $\Gamma (W, \0_W)  \to \Gamma (X_i, \mathcal{O}_{X_i})$, which is an isomorphism
since $W$ is normal. So $W \cong \mathrm{Spec}(\Gamma (X_i, 
\mathcal{O}_{X_i})), i=1,2$.
Therefore, the two resolutions $\pi_1$ and $\pi_2$ are equivalent.
\end{proof}
\begin{Cor}
Let $\0$ be a nilpotent orbit in $\mathfrak{sl}(n)$ with Jordan decomposition type $[d_1, \cdots, d_k]$. Suppose that $d_1 =2$, then any two
symplectic resolutions for $\overline{\0}$ are equivalent.
\end{Cor}
\begin{proof}
Since $d_1 = 2$, the dual partition of  $[d_1, \cdots, d_k]$ consists of two parts $[n-t,t]$, where $t= \#\{i|d_i =2\}$. 
$\overline{\0}$ has two symplectic resolutions, which are given by cotangent spaces of Grassmanians: $T^*Gr(n,t) \to \overline{\0}$ 
and $T^*Gr(n, n-t) \to \overline{\0}$, thus they are equivalent.
\end{proof}

Some interesting questions relating to derived categories for the two symplectic resolutions $T^*Gr(n,t) \to \overline{\0}$ 
and $T^*Gr(n, n-t) \to \overline{\0}$
are discussed in \cite{Na2}.

\begin{Exam} Here we give an example where a nilpotent orbit admits two non-isomorphic symplectic resolutions.
 Let $n \geq 2$ be an integer.  Consider the symplectic resolution $T^*\pit^n \xrightarrow{\pi} \overline{\0}_{min}$, where $\0_{min} = \0_{[2,1^n]}$
is the minimal nilpotent orbit in $\mathfrak{sl}(n+1, \cit).$  Now we perform a Mukai flop along the zero section $P \simeq \pit^n$ of 
$T^*\pit^n$, i.e. we first blow up $ T^*\pit^n  $ along $P$, then blow down along another direction to get another symplectic resolution 
$T^*\pit^n \xrightarrow{\pi^+} \overline{\0}_{min}. $ Notice that the birational map $(\pi^+)^{-1} \circ \pi: T^*\pit^n --\to T^*\pit^n  $
is not defined at the points of $P$. So $\pi$ and $\pi^+$ are not isomorphic.  In fact, the two symplectic resolutions come from non-conjugate
parabolic subgroups in $G$, one is the stabilizer of a line in $\cit^{n+1}$ and the other is the stabilizer of 
a vector subspace of codimension 1 in $\cit^{n+1}$.
\end{Exam}

\begin{Exam} \label{flag}
Here we give an example of a nilpotent orbit closure which admits three non-equivalent symplectic resolutions.
Let $\mathcal{O}$ be the nilpotent orbit 
of $\mathfrak{sl}(6)$ with Jordan decomposition type 
$[3,2,1]$. Then there are six non-conjugate polarizations 
$P_{(\sigma (1), \sigma (2), \sigma (3))}$ of 
$\mathcal{O}$, where $\sigma$ 
is a permutation of $\{1,2,3\}$. 
There are six  
non-isomorphic symplectic resolutions of 
$\bar{\mathcal{O}}$ corresponding to the 
six polarizations. Among these,  
the following pairs are equivalent resolutions by Proposition \ref{sl}: 
$$(T^*F(6,3,1), T^*F(6,5,3))$$  
$$(T^*F(6,3,2)), T^*F(6,4,3))$$ 
$$(T^*F(6,5,2)), T^*F(6,4,1)).$$ 
We now show that there are exactly three 
non-equivalent resolutions. Assume that 
two of three cotangent bundles $T^*F(6,3,1)$, 
$T^*F(6,3,2)$ and $T^*F(6,5,2)$ are equivalent 
resolutions of $\bar{\mathcal{O}}$.  
Let us consider the fibers of each resolution. 
Since the fibers with $\dim = 1/2 \dim T^*F$ are 
central fibers, if two resolutions are equivalent, 
then the corresponding flag manifolds are mutually 
isomorphic.    
We shall prove that this is absurd. 
We observe 
ample cones of these varieties. Since these varieties 
have Picard number two, they have at most 
two different fibrations. 
$F(6,3,1)$ has two fibrations 
$F(6,3,1) \to F(6,3)$ and $F(6,3,1) \to F(6,1)$. 
The first one is a $\pit^2$-bundle and 
the second one is a $Gr(5,2)$-bundle. 
$F(6,3,2)$ has two 
fibrations $F(6,3,2) \to F(6,3)$ and 
$F(6,3,2) \to F(6,2)$. The first one is a $\pit^2$-bundle 
and the second one is a $\pit^3$-bundle. 
$F(6,5,2)$ has two 
fibrations $F(6,5,2) \to F(6,5)$ and 
$F(6,5,2) \to F(6,2)$. The first one is a $Gr(5,2)$-bundle 
and the second one is a $\pit^3$-bundle. 
If two of these varieties are isomorphic, they should 
have three different fibrations, which is absurd. 

By Lemma \ref{equiv}, we see that neither two of $T^*F(6,3,1), T^*F(6,3,2)$ and $T^*F(6,5,2)$ are isomorphic.
\end{Exam}     

\section{Finiteness of symplectic resolutions and deformations}
We propose the following conjecture: 
\begin{Conj}
Let $W$ be a normal symplectic singularity. Then $W$ admits at most finitely many non-isomorphic symplectic resolutions.
\end{Conj}

Note that for nilpotent orbits, this conjecture is proved in Theorem \ref{nilthm}.
Here we prove the conjecture in the case of $\dim(W) = 4$.
\begin{Thm}
There are only finitely many non-isomorphic symplectic resolutions of a symplectic singularity $W$ of dimension 4. 
\end{Thm}
\begin{proof}
Fix a symplectic resolution $f: X \to W$.  Let $f^+: X^+ \to W$ be another symplectic resolution. 
Then, $X$ and $X^+$ are connected by a finite sequence of Mukai flops over $W$ 
(the existence of the flops follows from \cite{CMS} or \cite{WW}, and the termination of the flop sequence  follows from \cite{Mat}). 
Then we can apply the  argument of \cite{KM} to prove our theorem.     
\end{proof}

\begin{Exam} Let $A$ be an abelian surface and $\sigma: A \to A$ the involution $x \mapsto -x$. Then $A_0: = A/<\sigma>$ has 16 double points. Let
$B \to A_0$ be the minimal resolution. Then $\pi: B \times B \to A_0 \times A_0$ is a symplectic resolution.
 Notice that the 2-dimensional $\pi$-exceptional  fibers are isomorphic to $\pit^1 \times \pit^1$, thus no Mukai flop can be performed.
Thus $\pi$ is the unique symplectic resolution for $A_0 \times A_0$, up to isomorphisms.  \end{Exam}

\begin{Exam} Let $f: \mathrm{Hilb}^2(S) \to \mathrm{Sym^2(\bar{S})}$ be the symplectic resolution considered in Example \ref{exam}.
The only 2-dimensional $f$-exceptional fiber is $f^{-1}(0) = \pit^2 \cup F$. We can perform only once Mukai flop to $f$, thus
$\mathrm{Sym^2(\bar{S})}$ admits exactly two non-isomorphic symplectic resolutions: $f$ and $\bar{u} \circ f$. \end{Exam}

Recall that a deformation of a variety $X$ is a 
 flat morphism 
$\mathcal{X} \xrightarrow{p} S$ from a variety $\mathcal{X}$ to a 
pointed smooth connected curve $0 \in S$ such that $p^{-1}(0) \cong 
X$. Moreover, a deformation of a proper morphism $f: X \to Y$ 
is a proper $S$-morphism $F: \mathcal{X} \to \mathcal{Y}$, where 
$\mathcal{X} \to S$ is a deformation of $X$ and 
$\mathcal{Y} \to S$ is a deformation of $Y$. 
 
Two varieties $X_1$ and $X_2$ are said {\em deformation equivalent} if there is a flat morphism $\mathcal{X} \xrightarrow{p} S$ from a variety 
$\mathcal{X}$ to a connected (not necessarily irreducible) 
curve $S$ such that  
 $X_1$ and $X_2$ are isomorphic to two fibers of $p$.  As to the relation between two symplectic resolutions, we have the following:

\begin{Conj}\label{Conj2}
Let $W$ be a normal symplectic singularity. Then for any two symplectic resolutions $f_i: X_i \to W, i=1, 2$, 
there are deformations $\mathcal{X}_i \stackrel{F_i}\to 
\mathcal{W}$ of $f_i$ such that, for $s \in S\setminus 0$, $F_{i,s}: 
\mathcal{X}_{i,s} \to \mathcal{W}_s$ are isomorphisms. In particular, 
$X_1$ and $X_2$ are deformation equivalent. 
\end{Conj}

If $W$ is a projective symplectic variety (with singularities), 
then we have Kuranishi spaces $\mathrm{Def}(W)$ and 
$\mathrm{Def}(X_i)$ for $W$ and $X_i$. Since $W$ has only 
rational singularities, we have the maps $(f_i)_*: 
\mathrm{Def}(X_i) \to \mathrm{Def}(W)$. 
By Theorem 2.2, \cite{Nam}, the Kuranishi spaces are all 
non-singular and $(f_i)_*$ are finite coverings. Now take a 
map $\Delta \to \mathrm{Def}(W)$ from a 1-dimensional disk 
such that this map factors through both $\mathrm{Def}(X_i)$. 
By pulling back the semi-universal families by this map, 
we have three flat families of varieties. If we take the map 
sufficiently general, then these families give the desired ones 
in the conjecture. One can say more. D. Huybrechts in \cite{Huy} 
proved that if two compact hyper-K\"ahler manifolds $X_1$ and $X_2$ 
are birationally equivalent, then they are deformation equivalent. 
Here we do not need the intermediate variety $W$ any more. 

Let us return to our local case. When $W$ is an isolated singularity, 
we also have the Kuranishi spaces for $W$ and $X_i$. Moreover, 
by \cite{CMS} and \cite{WW}, $f_i$'s give a Mukai flop in this case. 
Then one can show the Conjecture applying the deformation theory 
as well as the projective case. The problem is when $W$ is 
not an isolated singularity. We do not have appropriate spaces 
like the Kuranishi spaces any more. Sometimes, the formal approach 
could be possible, but its convergence is a difficult problem.  
 D. Kaledin proved this conjecture under some hypothesis
in \cite{Kal}. 
For the last statement of the conjecture, we proved in \cite{Fu2} 
that $X_1$ is deformation equivalent to $X_2$ when they are  symplectic resolutions of nilpotent orbit 
closures in a classical simple complex Lie algebra.
Here we prove Conjecture \ref{Conj2} for nilpotent orbit closures in $\g = \mathfrak{sl}(n)$. The construction is elementary and it may
be of independent interest (see \cite{Na2}). 

\begin{Thm}
Let $\0$ be a nilpotent orbit in $\mathfrak{sl}(n)$ with Jordan decomposition type
$[d_1, \cdots, d_k]$. Then  Conjecture \ref{Conj2} holds for $\overline{\0}$.
\end{Thm}
\begin{proof}
Let $[s_1, \cdots, s_m]$ be the dual partition of $[d_1, \cdots, d_k]$. 
The polarizations of $\0$ are $P_{s_{\sigma(1)}, \cdots, s_{\sigma(m)}}, \sigma \in \Sigma_m.$
Define $F_\sigma : = SL(n)/P_{s_{\sigma(1)}, \cdots, s_{\sigma(m)}}. $ 
Let $$ \tau_1  \subset \cdots \subset  \tau_{m-1} \subset
\cit^n \otimes_{\cit}\mathcal{O}_{F_{\sigma}}$$ 
be the universal subbundles on $F_{\sigma}$. 
A point of $T^*F_\sigma$ is expressed as a pair 
$(p, \phi)$ of $p \in F_{\sigma}$ and $\phi \in 
\mathrm{End}(\cit^n)$ such that 
$$ \phi(\cit^n) \subset \tau_{m-1}(p), \cdots,  
\phi(\tau_2(p)) \subset 
\tau_1(p), \phi(\tau_1(p)) = 0.$$ 
The Springer resolution 
$$ s_{\sigma}: T^*F_{\sigma} \to \bar{\mathcal{O}}$$ 
is defined as $s_{\sigma}((p, \phi)) := \phi$.
 
First, we shall define a vector bundle 
$\mathcal{E}_{\sigma}$ over $F_{\sigma}$ and an exact sequence 
$$ 0 \to T^*F_{\sigma} \to \mathcal{E}_{\sigma} 
\stackrel{\eta_{\sigma}}\to 
\mathcal{O}_{F_{\sigma}}^{\oplus m-1} \to 0. $$ 
Let $T^*F_{\sigma}(p)$ be the cotangent space of 
$F_{\sigma}$ at $p \in F_{\sigma}$. Then, for a suitable 
basis of $\cit^n$, $T^*F_{\sigma}(p)$ consists of the matrices 
of the following form 
$$ 
\begin{pmatrix} 
0 & * & \cdots &* \\
0 & 0 &  \cdots & *\\
\cdots & & & \cdots \\
0 & 0 &  \cdots& 0 \\ 
\end{pmatrix}.
$$ 
Let $\mathcal{E}_{\sigma}(p)$ be the vector 
subspace of $\mathfrak{sl}(n)$ consisting of 
the matrices $A$ of the following 
form
$$ 
\begin{pmatrix} 
a_{\sigma (1)} & * & \cdots & * \\
0 & a_{\sigma (2)} &\cdots  & *\\
\cdots & & & \cdots \\
0 & 0 & \cdots  &a_{\sigma (m)} \\ 
\end{pmatrix},
$$
where $a_i := a_i I_{s_i}$ and $I_{s_i}$ is 
the identity matrix of 
the size $s_i \times s_i$. 
Since $A \in \mathfrak{sl}(n)$, 
$\Sigma_i s_i a_i = 0$. 
We define a map $\eta_{\sigma}(p): \mathcal{E}_{\sigma}(p) 
\to {\cit}^{\oplus m-1}$ as 
$\eta_{\sigma}(p)(A) := (a_1, a_2, \cdots, a_{m-1}).$ 
Then we have an exact sequence of 
vector spaces 
$$ 0 \to T^*F_{\sigma}(p) \to \mathcal{E}_{\sigma}(p) 
\stackrel{\eta_{\sigma}(p)}\to {\cit}^{\oplus 
m-1} \to 0. $$  
We put $\mathcal{E}_{\sigma} := \cup_{p \in F_{\sigma}}
\mathcal{E}_{\sigma}(p)$. Then $\mathcal{E}_{\sigma}$ 
becomes a vector bundle over $F_{\sigma}$, and 
we get the desired exact sequence. 
Note that we have a morphism 
$$\eta_{\sigma}: 
\mathcal{E}_{\sigma} \to \cit^{m-1}. $$   
 
Next, let $\overline{N} \subset \mathfrak{sl}(n)$ be the set of 
all matrices which is conjugate to a matrice
of the following form:
$$ 
\begin{pmatrix} 
b_{1} & * & \cdots & * \\
0 & b_{2} &\cdots  & *\\
\cdots & & & \cdots \\
0 & 0 & \cdots  & b_{m} \\ 
\end{pmatrix},
$$
where $b_i = b_i I_{s_i}$ and  
$I_{s_i}$ is the identity matrix of order $s_{i}$. Furthermore 
the zero trace 
condition requires $\sum_i s_i b_i =0$. 
For $A \in \overline{N}$, let $\phi_A(x) := \mathrm{det}
(xI - A)$ be the characteristic polynomial of 
$A$. Let $\phi_i(A)$ be the coefficient of 
$x^{n-i}$ in $\phi (A)$. 
We define the characteristic map $ch: \overline{N} 
\to \cit^{n-1}$ by $ch(A) := (\phi_2(A), ..., 
\phi_n(A))$. Note that $\phi_1(A) = 0$. 

Let us consider the vector 
$\mathbf{a} = (a_1, a_2, ...., a_m)$ 
of length $n$ where $a_i$ appear exactly $s_i$ times.  
Define $\phi_{i,\mathbf{a}}$ to be the $\phi_i(A)$ 
for the following diagonal matrix $A$ of the size $n \times n$  
$$ 
\begin{pmatrix} 
a_1 & 0 & \cdots & 0 \\
0 & a_2 & 0 & \cdots \\
\cdots & & & \cdots \\
0 & 0 & \cdots  & a_m \\ 
\end{pmatrix},
$$
where $a_{i} = a_{i}I_{s_i}$ with $s_i \times s_i$ 
identity matrix $I_{s_i}$.  
For $(a_1, a_2, ..., a_{m-1}) \in \cit^{m-1}$ 
we put $\mathbf{a}' := (a_1, ..., a_{m-1}, 
-\Sigma_{i=1}^{m-1} s_i a_i/s_m)$, and we  
define a map 
$$\pi: \cit^{m-1} \to \cit^{n-1}$$ 
by $\pi(a_1, ..., a_{m-1}) = 
(\phi_{2,\mathbf{a}'}, ..., \phi_{n, \mathbf{a}'})$.   
Pulling back  
$ch: \overline{N} \to \cit^{n-1}$ by $\pi$, we have 
$$ ch': \overline{N}' \to \cit^{m-1}. $$    

Each point of $\mathcal{E}_{\sigma}$ is 
expressed as a pair of a point $p \in 
F_{\sigma}$ and $\phi \in \mathrm{End}(\cit^n)$. 
Now we define 
$$ \bar{s}_{\sigma}: \mathcal{E}_{\sigma}
 \to \overline{N} $$ 
as $\bar{s}_{\sigma}(p, \phi) = \phi$. 
This map is a generically finite morphism. 
Since $ch \circ \bar{s}_{\sigma} = 
\pi \circ \eta_{\sigma}$, we have a morphism 
$$ \bar{s}'_{\sigma}: 
\mathcal{E}_{\sigma} \to \overline{N}'. $$ 
Let $\widetilde{N}$ be the normalization of 
$\overline{N}'$ and let $f$ be the 
composite: $\widetilde{N} \to 
\overline{N}' \stackrel{ch'}\to \cit^{m-1}$.
Then $\bar{s}'_{\sigma}$ 
factors through $\widetilde{N}$ and we have a 
morphism 
$$ \tilde{s}_{\sigma}: \mathcal{E}_{\sigma} 
\to \widetilde{N}. $$ 
Now, $\tilde{s}_{\sigma}$ becomes a birational 
morphism. Moreover, for a general point 
$\mathbf{t} \in \cit^{m-1}$, 
$\tilde{s}_{\sigma, \mathbf{t}}$ is an isomorphism. 
The flat deformations 
$$  \mathcal{E}_{\sigma} \stackrel{\tilde{s}_{\sigma}}
\to \widetilde{N} \stackrel{f}\to \cit^{m-1}$$ 
give desired deformations in the conjecture.
\end{proof}

\section{An example}

In this section we construct a symplectic singularity $W$ of dim 4 which has two non-equivalent 
symplectic resolutions. 
We already have such examples by the nilpotent 
orbit construction (cf. \ref{flag}). 
But here we introduce another construction. 
Our construction is elementary.

A similar example has also been constructed by J. Wierzba (Section 7.2.3 \cite{Wi2}), using a
different approach. Finally we note that such an example 
can be constructed by hyper-K\"ahler quotients 
\cite{Got}.   

\subsection{The idea}\label{idea}  Let $f: V \to W$ 
be a symplectic resolution such that:

(i). for some point 
$0 \in W$, $f^{-1}(0) = \pit^2 
\cup \Sigma_1$, where $\pit^2 
\cap \Sigma_1$ is a line on $\pit^2$ 
and, is, at the same time, a negative 
section of $\Sigma_1$;

(ii). the singular locus 
$\Sigma$ of $W$ is 2-dimensional. And for $p 
\in \Sigma$ such that $p \ne 0$, 
$(W, p) \cong (A_1-\mathrm{surface\ singularity}) 
\times (\cit^2, 0)$. 

Over such a point 
$p$, $f$ will become the minimal resolution.  
Now flop $V$ along $\pit^2$; then we get 
a new symplectic resolution $f^+: V^+ \to W$ 
such that ${f^+}^{-1}(0) = \pit^2 
\cup \pit^2$ where two $\pit^2$ 
intersect in one point. Then, it is clear 
that the two symplectic resolutions are not 
equivalent. In fact, if they are equivalent, 
then there should be an isomorphism 
$V \cong V^+$ which sends $f^{-1}(0)$ 
isomorphically onto ${f^+}^{-1}(0)$. But 
this is absurd.  

\subsection{ Construction of the example}

\subsubsection{Set-up} Let $\bar{S}$ be the germ of an $A_2$-surface 
singularity and let $\pi: S \to \bar{S}$ be 
its minimal resolution with exceptional curves 
$C$ and $D$. There are 
natural birational morphisms 
$$ \mathrm{Hilb}^2(S) \xrightarrow{\nu} 
\mathrm{Sym}^2(S) \to \mathrm{Sym}^2(\bar{S}).$$ 
We denote by $g: \mathrm{Hilb}^2(S) \to 
\mathrm{Sym}^2(\bar{S})$ the composition. 
$\mathrm{Sym}^2(S)$ contains 
$\mathrm{Sym}^2(C)$ and $\mathrm{Sym}^2(D)$. 
Let $P_C$ and $P_D$ be their proper transforms 
on $\mathrm{Hilb}^2(S)$. Note they are 
isomorphic to $\pit^2$. 
Let us consider the double cover 
$\alpha: S \times S \to \mathrm{Sym}^2(S)$. 
Then $\alpha (C \times D) = \alpha (D \times C)$ 
and $\alpha (C \times D) \cong C \times D$. 
Let $Q$ be the proper transform of 
$\alpha (C \times D)$ on $\mathrm{Hilb}^2(S)$. 
Now $Q$ is isomorphic to the one point blow-up 
of $C \times D$. If $p := C \cap D$ in $S$, then 
the center of the blowing-up is $(p, p) \in 
C \times D$. Let $l_C \subset Q$ (resp. $l_D 
\subset Q$) be 
the proper transform of $C \times \{p\}$ 
(resp. $\{p\} \times D$) by the blowing-up. 
Let $e \subset Q$ be the exceptional curve. 
Then $l_C$, $l_D$ and $e$ are $(-1)$-curves 
of $Q$ with $(l_C, l_D) = 0$, $(l_C, e) 
= (l_D, e) = 1$.   
The relationship between $P_C$, 
$P_D$ and $Q$ are the following. 
\vspace{0.12cm}

(i) $P_C$ and $P_D$ are disjoint. 

(ii) $Q$ intersects both $P_C$ and $P_D$. 
 
(iii) In $Q$, $Q \cap P_C$ coincides with $l_C$ 
and $Q \cap P_D$ coincides with $l_D$. 

(iv) In $P_C$, $Q \cap P_C$ is a line, and, in 
$P_D$, $Q \cap P_D$ is a line.  
\vspace{0.12cm}

Let $E \subset \mathrm{Hilb}^2(S)$ be the 
exceptional divisor of the birational 
morphism $\nu : \mathrm{Hilb}^2(S) 
\to \mathrm{Sym}^2(S)$. Let $E_C := E \cap 
\nu^{-1}(\mathrm{Sym}^2(C))$ and 
$E_D := E \cap \nu^{-1}(\mathrm{Sym}^2(D))$. 
$E_C$ is a $\pit^1$-bundle over 
the diagonal $\Delta_C \subset \mathrm{Sym}^2(C)$. 
Let $f_C$ be a fiber of this bundle. 
Similarly, $E_D$ is a $\pit^1$-bundle 
over $\Delta_D \subset \mathrm{Sym}^2(D)$, and 
let $f_D$ be its fiber.  
Note that 
$$g^{-1}(0) = Q \cup P_C \cup P_D 
\cup E_C \cup E_D. $$

\subsubsection{Mukai flop} Flop $\mathrm{Hilb}^2(S)$ along the 
center $P_C$ 
to get a new 4-fold $V$. We denote by 
$P'_C \subset V$ the center of this flop. 
There is a birational morphism 
$g^+: V \to \mathrm{Sym}^2(\bar{S})$. 

Let $P'_D \subset V$ 
be the proper transform of $P_D$, and let 
$Q' \subset V$ be the proper transform of $Q$. 
Since $P_D$ is disjoint from $P_C$, $P'_D$ 
is naturally isomorphic to $P_D$; hence 
$P'_D \cong \pit^2$. 
On the other hand, $Q'$ is isomorphic to 
the blowing down of $Q$ along $l_C$. 
Now $Q'$ becomes the Hirzebruch surface 
$\Sigma_1$. The intersection $P'_D \cap Q'$ 
is a line of $P'_D$, and is a negative 
section of $Q' \cong \Sigma_1$. 
On the other hand, $Q'$ and $P'_C$ intersect 
in one point. 
Let $E'_C \subset V$ (resp. $E'_D \subset V$) be 
the proper transform of $E_C$ (resp. $E_D$). 

\subsubsection{Idea} We shall construct a birational contraction 
map $f: V \to W$ over $\mathrm{Sym}^2(\bar{S})$ 
such that, in $(g^+)^{-1}(0)$, 
$P'_D$ and $Q'$ are contracted 
to a point by $f$, $E'_C$ is contracted along 
the ruling to a curve, and both  
$E'_D$ and $P'_C$ are birationally mapped onto their 
images. We put $f(P'_D \cap Q') = q$ and let 
$W^0$ be a sufficiently small open neighborhood 
of $q \in W$, and let $V^0 := f^{-1}(W^0)$. 
Then $f^0 (:= f\vert_{V^0}): V^0 \to W^0$ satisfies  
the conditions 
of section \ref{idea}. Let $(f^0)^+ : (V^0)^+ \to W^0$ be 
another symplectic 
resolution obtained by flopping $P'_D$. Then 
$f^0$ and $(f^0)^+$ are not equivalent.

\subsubsection{The construction of $f$}
Let $\mu: B(S \times S) \to S \times S$ 
be the blowing-up along the diagonal $\Delta_S$. 
Let $F$ be the exceptional divisor of 
the blowing-up. We have a double cover 
$\tilde{\alpha} : B(S \times S) \to 
\mathrm{Hilb}^2(S)$. We can write 
$$\tilde{\alpha}_*\mathcal{O}_{B(S\times S)} 
= \mathcal{O}_{\mathrm{Hilb}^2(S)} \oplus 
M$$ for some $M \in \mathrm{Pic}(\mathrm{Hilb}^2
(S))$. Note that $M^{\otimes 2} = \mathcal{O}(-E)$. 
Choose $L \in \mathrm{Pic}(S)$ in such a way that 
$(L.C) = 0$ and $(L.D) = 1$. The line bundle 
$\mu^*(p_1^*L \otimes p_2^*L)$ on 
$B(S \times S)$ can be written as the pull-back by 
$\tilde{\alpha}$ 
of a line bundle $N$ on $\mathrm{Hilb}^2(S)$. 
Define 
$$ \mathcal{L} := N \otimes M. $$ 
Then we have  
$$ (\mathcal{L}. e) = 1, (\mathcal{L}. l_C) = -1, 
(\mathcal{L}. l_D) = 0 $$ 
$$ (\mathcal{L}. f_C) = 1, (\mathcal{L}. f_D) = 1.$$ 

We have the following situation after the flop along 
$P_C$.  \vspace{0.12cm}
 
(i) The proper transform $e'$ of $e$ 
is a ruling of $Q' \cong \Sigma_1$.

(ii) The proper transform $l'_D$ of $l_D$ 
is a negative section of $Q' \cong \Sigma_1$, 
and at the same time, is a line of $P'_D$. 
\vspace{0.12cm}

Let $f'_C$ be the proper transform of $f_C$, and 
let $f'_D$ be the proper transform of $f_D$. Then, 
for the proper transform $\mathcal{L}' \in 
\mathrm{Pic}(V)$ of $\mathcal{L}$, we have 
$$ (\mathcal{L}'. e') = 0, (\mathcal{L}'. l'_D) 
= 0, (\mathcal{L}'. f'_C) = 0, 
(\mathcal{L}'. f'_D) = 1. $$ 
Moreover, for a line $l$ of $P'_C$, we 
see that $(\mathcal{L}'.l) = 1$ because 
$\mathrm{Hilb}^2(S) -- \to V$ is the 
flop along $P_C$ and $(\mathcal{L}. l_C) 
= -1$.  
These implies that $\mathcal{L}'$ is 
$g^+$-nef (and, of course, $g^+$-big).  
Since $g^+$ is a crepant resolution, by 
the base point free theorem, 
$\mathcal{L}'^{\otimes n}$ is $g^+$-free 
for a sufficiently large $n$. By this line 
bundle we define $f: V \to W$. 
An irreducible curve on $V$ is contracted to 
a point if and only if it has no intersection 
number with $\mathcal{L}'$. Since 
$(\mathcal{L}'. l'_D) = 0$, $f$ contracts $P'_D$  
to a point by (ii). Moreover, since 
$(\mathcal{L}'. e') = 0$, $f$ contracts $Q'$ to 
the same point. Finally, since $(\mathcal{L}'.f'_C) 
= 0$, $f$ contracts every ruling of $E'_C$ to 
ponits. Similarly we can check that $P'_C$ and 
$E'_D$ are birationally mapped onto their images 
by $f$.   

\subsubsection{Detailed description of $f$} 
Among the irreducible components of 
${g^+}^{-1}(0)$, $Q'$, $P'_D$ and $E'_C$ 
are $f$-exceptional. 
The birational morphism $f$ factorize 
$g^+$ as 
$$ V \stackrel{f}\to W \stackrel{h}\to 
\mathrm{Sym}^2(\bar{S}). $$ 
Then $h^{-1}(0)$   
consists of two components; one of them is 
$f(P'_C) \cong \pit^2$ and another one 
is $f(E'_D)$, which is the blow-down of 
$E'_D \cong \Sigma_4$ along the negative 
section. These two components intersect in 
one point. Note that $f(E'_C)$ is a conic 
on $f(P'_C)$.    
\begin{Rque}
It follows from Lemma \ref{equiv} that for the two symplectic resolutions $V \to W$ and $V^+ \to W$, $V$ is not isomorphic to $V^+$.
\end{Rque}

\quad \\
\quad\\
Baohua Fu\\
Labortoire J. A. Dieudonn\'e, Parc Valrose \\ 06108 Nice cedex 02, FRANCE \\
baohua.fu@polytechnique.org

\quad\\
Yoshinori Namikawa \\
Departement of Mathematics, Graduate School of Science, Kyoto University,
Kita-shirakawa Oiwake-cho, Kyoto, 606-8502, JAPAN \\
namikawa@kusm.kyoto-u.ac.jp

\end{document}